# Space Exploration Architecture and Design Framework for Commercialization

Hao Chen[*]
*Georgia Institute of Technology, Atlanta, GA, 30332*

Melkior Ornik[†,‡]
*University of Illinois Urbana-Champaign, 104 S. Wright St, Urbana, IL 61801*

and
Koki Ho[§]
*Georgia Institute of Technology, Atlanta, GA, 30332*

**The trend of space commercialization is changing the decision-making process for future space exploration architectures, and there is a growing need for a new decision-making framework that explicitly considers the interactions between the mission coordinator (i.e., government) and the commercial players. In response to this challenge, this paper develops a framework for space exploration and logistics decision making that considers the incentive mechanism to stimulate commercial participation in future space infrastructure development and deployment. By extending the state-of-the-art space logistics design formulations from the game-theoretic perspective, the relationship between the mission coordinator and commercial players is first analyzed, and then the formulation for the optimal architecture design and incentive mechanism in three different scenarios is derived. To demonstrate and evaluate the effectiveness of the proposed framework, a case study on lunar habitat infrastructure design and deployment is conducted. Results show how total mission demands and in-situ resource utilization (ISRU) system performances after deployment may impact the cooperation among stakeholders. As an outcome of this study, an incentive-based decision-making framework that can benefit both the mission coordinator and the commercial players from commercialization**

---

[*] Postdoctoral Fellow, Daniel Guggenheim School of Aerospace Engineering, AIAA Student Member.
[†] Assistant Professor, Department of Aerospace Engineering.
[‡] Assistant Professor, Coordinated Science Laboratory.
[§] Assistant Professor, Daniel Guggenheim School of Aerospace Engineering, AIAA Member.

**is derived, leading to a mutually beneficial space exploration between the government and the industry.**

## Nomenclature

| | | |
|---|---|---|
| $\mathcal{A}$ | = | set of arcs |
| $\mathcal{C}$ | = | set of commodities |
| $\mathcal{C}_C$ | = | set of continuous commodities |
| $\mathcal{C}_D$ | = | set of discrete commodities |
| $c$ | = | cost coefficient |
| $d$ | = | mission demand |
| $G$ | = | commodity transformation matrix |
| $H$ | = | concurrency constraint matrix |
| $I_{sp}$ | = | specific impulse, s |
| $\mathcal{J}$ | = | space mission cost |
| $\mathcal{K}$ | = | set of commercial players |
| $l$ | = | concurrency constraint index |
| $\mathcal{N}$ | = | set of nodes |
| $N$ | = | total number of players |
| $Q$ | = | baseline mission cost |
| $r$ | = | disagreement point utility |
| $\mathcal{T}$ | = | set of time steps |
| $U$ | = | total utility |
| $u$ | = | mission utility |
| $\mathcal{V}$ | = | set of spacecraft |
| $W$ | = | set of time windows |
| $x$ | = | commodity variable |
| $\mathcal{Z}$ | = | set of game players |
| $\Delta t$ | = | time of flight, day |

| $\Delta V$ | = | change of velocity, km/s |
| $\alpha$ | = | participation coefficient |
| $\theta$ | = | incentive coefficient |
| $\Omega$ | = | feasible domain of the participation and incentive coefficients |

*Subscripts*

| $i$ | = | node index |
| $j$ | = | node index |
| $k$ | = | commercial player index |
| $m$ | = | commercial player index |
| $t$ | = | time step index |
| $v$ | = | spacecraft index |
| $z$ | = | game player index |

## I. Introduction

THERE is a growing trend of space commercialization; more government and commercial entities exhibit their interest in participating in future large-scale space exploration. Each entity has potentially different mission objective preferences and technology advantages, and the coordination among these entities will become a critical component of space exploration and design practice. This trend is also closely related to the needs of campaign-level mission design for future deep space explorations beyond Earth orbits, where space logistics infrastructures, such as in-situ resource utilization (ISRU) systems [1], will play a critical role in reducing space mission cost leveraging mission interdependencies.

Campaign-level space exploration and logistics design methodologies have been extensively studied recently. Earlier works focused on discrete-event simulation environments [2], followed by the graph-theoretic or network-based modeling of space logistics [3-5]. Based on the space logistics network, dynamic space mission design optimization methodologies have been established using the time-expanded generalized multi-commodity network flow model [6-8]. Depending on the mission designers' needs, the campaign-level space logistics model has been extended to the concurrent optimization of the mission scheduling and the optimal allocation of vehicles (high-thrust or low-thrust) to each mission leg [9] as well as the optimal subsystem-level design of the deployed infrastructure such as ISRU plants [10]. Furthermore, for a long-horizon campaign with regularly repeating demands, a partially periodic time-expanded network has been developed to simultaneously optimize the infrastructure deployment phase and operational phase [11].

Despite the rich literature in the space logistics field, all the aforementioned space logistics mission design methods have assumed that a single stakeholder has top-down control over all systems and resources during the space campaign. The optimization formulation proposed in these studies can lead to a theoretically optimal mission design, which would be helpful if one entity (e.g., NASA) controls the entire process. However, these theoretically optimal solutions are not necessarily practical because they do not consider the interaction between the coordinator and commercial players. While some studies have looked into particular commercial opportunities in different industrial fields for space exploration, such as the satellite industry [12,13] and the commercialization of low-Earth orbit (LEO) [14], they did not take into account space architecture and design as part of the trade space, and therefore cannot be used for mission architecture decision making for commercialization. Other studies have also looked into the commercial suitability of space infrastructure design from the modularity perspective [15], but they only focused on limited commercialization perspectives and did not explicitly consider the incentive and objectives for commercial players.

A related topic studied in the literature is the federated system (i.e., system-of-systems, SoS), which has been used as one of the important ways to stress cooperation and competition among multiple actors. The idea of a federated system has been proposed and applied to federated satellite systems [16], the government extended enterprises SoS [17], and the multi-actor space architecture commercialization [18]. However, these studies solely focused on allocating resources to each player rather than considering the architecture design from the commercialization perspective. An effective federation needs to allow players to actively decide to participate in or leave the federation considering their mission objective, available resources, architecture design, and potential benefits. This scenario would require the government/coordinator to offer incentives and perform architecture design from the commercialization perspective, which is the focus of this paper.

The problem of interest is how to perform effective architecture design together with the incentive mechanism to achieve their goal with minimum resources. To solve this problem, this paper proposes an architecture and design decision-making framework that can incentivize commercial players to join the space exploration enterprise. The proposed architecture and design decision-making method aims to build a mutually beneficial relationship between the government mission coordinator (e.g., NASA) and the industry. Incentive mechanism design has been widely studied in various applications, including mobile phone participatory sensing [19-21], corporate entrepreneurship [22], and health risk assessment [23]. Specifically, a recent study focusing on infrastructure deployment [24] established an incentive design method for the mobile network market. Leveraging the approach in the incentive mechanism design field and extending the state-of-the-art space logistics and mission design approaches, we propose an architecture design framework for commercial participation in space infrastructure development and deployment. Building upon our previous work [25], we establish the framework based on the Nash bargaining solution and analyze its properties.

Depending on the mission planning circumstances, we discuss the different roles of the mission coordinator and the incentive design properties. There are two problems in the architecture design framework to solve. The first is the deployment task assignment problem to determine the mission demand for each player. The second is the incentive value optimization problem to distribute the total utility of the system through incentives. These two problems lead to three mission scenarios to be discussed in this paper: 1) when we can make decisions for both task assignment and incentive values (i.e., the global optimization can be solved); 2) when we are given a pre-determined deployment demand for each player (i.e., only the incentive value optimization problem is considered); 3) when we are given the fixed incentive budget allocation (i.e., only the task assignment problem is considered).

The major contributions of this paper will connect the incentive design for commercialization with the existing space mission design and logistics framework. It introduces a new perspective to analyze international and public-private federation with a diverse range of space activities and mission objectives from different space-faring countries and companies. We believe that the proposed framework will be a foundational work for future space exploration mission architecture and design through collaboration between government, international, and commercial players.

The remainder of this paper is organized as follows. Section II first introduces the game theory model settings and definitions. In Sec. III, we propose the architecture design framework based on game theory, particularly the Nash bargaining solution, and analyze its properties. Three different mission scenarios are discussed, and the corresponding solving approaches are established. The performance and comparison of the proposed architecture design framework are then demonstrated in Sec. IV through a lunar exploration campaign case study. Finally, Sec. V concludes the paper and discusses future works.

## II. Game Theory Model

In this section, we introduce the problem formulation and the game theory model. We consider a space logistics problem, where a government or commercial entity as a player has a mission demand to deploy some space infrastructure or deliver a certain amount of payload to a designated orbit or the lunar surface. We assume the total mass of the infrastructure or payload to be transported is $D$, in the unit of kg. This player can either complete this mission by itself only using its resources or relying on other players (i.e., other active government or commercial entities) to satisfy the mission demand. Other players may have deployed propellant depots, ISRU architectures, and lunar habitats that can complete the transportation mission at a lower cost. However, these commercial players also have their own mission goals. Transportation vehicle capacities and availabilities limit the available transportation capability for potentially extra commercial missions. We assume that all players are selfish but rational. Hence, they will only participate in the infrastructure deployment mission when it is beneficial. We define a player as a coordinator/planner who has a specific infrastructure deployment mission demand and would like to leverage other players' capability and resources to complete the space mission. A two-player incentive design example is shown in Fig. 1.

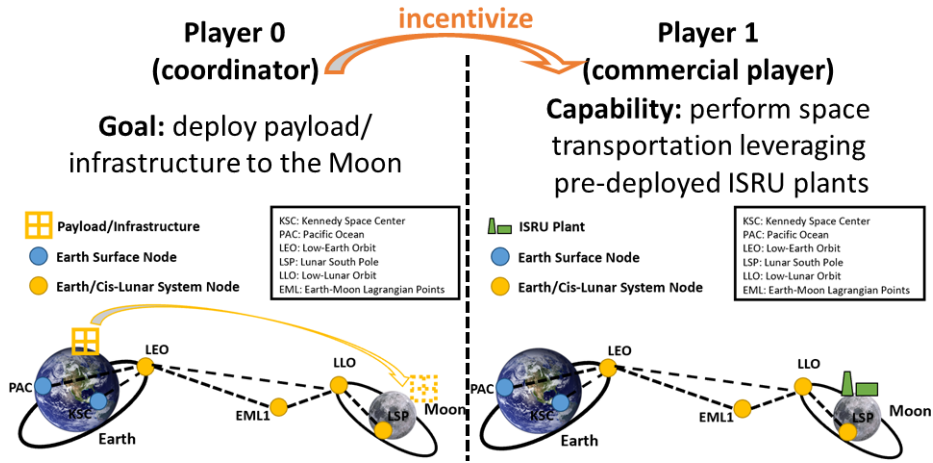

**Fig. 1 Incentive Design Problem Example**

For simplicity, in this paper, we consider a scenario with only one coordinator, labeled as player 0. All other players are commercial players, labeled as players 1, 2, etc. The problem is how to leverage the infrastructure owned by commercial players to support the space missions of the coordinator. For more general cases with multiple coordinators, we can always simplify the problem by considering the incentive design for each coordinator independently.

We define $\mathcal{K} = \{\text{commercial players}\}$ as the set of commercial players. Then, the total number of players $N$ in this game is the cardinality of $\mathcal{K}$ for commercial players plus one coordinator, $N = |\mathcal{K}| + 1$. According to available resources and mission demands, the coordinator can complete the transportation mission partially and complete the rest of the structure deployment relying on other players. We define a participation coefficient $\alpha_k \in [0,1]$ for each commercial player, denoting the fraction of mission demand to be deployed by the commercial player $k$. In this paper, mission demands are measured by the mass of payload in the unit of kg. Based on the above definition, the transportation mission demand for the coordinator is $(1 - \sum_{k \in \mathcal{K}} \alpha_k)D$; whereas the mission demand for the player $k$ is $\alpha_k D$. We assume that the total mission cost to complete the entire infrastructure deployment mission by the coordinator itself is $Q$, which we call it the baseline mission cost. The baseline cost to deliver 1 kg of infrastructure by the coordinator is $Q/D$. Then, we can define an incentive coefficient $\theta_k$. The feasible region of $\theta_k$ would need to satisfy certain constraints depending on $\alpha_k$; we denote the vectors of $\alpha_k$ and $\theta_k$ by $\boldsymbol{\alpha}$ and $\boldsymbol{\theta}$, and denote their feasible region as $\boldsymbol{\Omega}$ such that $(\boldsymbol{\alpha}, \boldsymbol{\theta}) \in \boldsymbol{\Omega}$; this domain $\boldsymbol{\Omega}$ will be discussed further later. As a rational coordinator, player 0 is only willing to pay the incentive no more than the baseline mission cost. The incentive to deliver 1 kg of infrastructure by the player $k$ can be denoted by $\theta_k Q/D$.

As a result, our decision variables include $\alpha_k$ and $\theta_k$. We define $u_o(\boldsymbol{\alpha}, \boldsymbol{\theta})$ as the utility of the coordinator. We also define $u_{p_k}(\alpha_k, \theta_k)$ as the utility of the commercial player $k$. The utility is defined as the mission cost savings for the coordinator or the profit for the commercial player. For the coordinator, the mission cost saving is calculated with respect to its baseline mission cost $Q$. Therefore, the utility function of the coordinator can be expressed as

$$u_o(\boldsymbol{\alpha}, \boldsymbol{\theta}) = Q - \mathcal{J}_o\left(1 - \sum_{k \in \mathcal{K}} \alpha_k\right) - \sum_{k \in \mathcal{K}} \alpha_k \theta_k Q \tag{1}$$

And the utility function of the commercial player $k$ can be written as

$$u_{p_k}(\alpha_k, \theta_k) = \alpha_k \theta_k Q - \mathcal{J}_{p_k}(\alpha_k) \tag{2}$$

where $\mathcal{J}_o$ and $\mathcal{J}_p$ are the space mission costs to complete the assigned portion of space infrastructure deployment for the coordinator and commercial players. Both are functions of only the participation coefficient $\alpha$ because the total infrastructure deployment demand $D$ is assumed as a known constant before the optimization. By definition, $\mathcal{J}_o(1) = Q$. These cost functions can be found numerically through space logistics optimization methods. Studies also have been done to establish analytical expressions for the mission cost functions [25, 26], but the accuracy of these methods depends on the actual mission scenarios. Note that the utility functions defined in Eqs. (1) and (2) have an implicit assumption that the incentive exchange between the coordinator and commercial players is frictionless. Friction losses (such as overhead) were not considered in this paper. Moreover, we also assume that the utilities of all players are fungible quantities like currency. More general utility functions of the coordinator and commercial players that cannot be exchanged among players are not considered in this framework.

We know that the utilities for all players are always non-negative; otherwise, they would not participate in the space mission. This leads to the constraints for the incentive coefficient $\boldsymbol{\theta}$:

$$Q - \mathcal{J}_o\left(1 - \sum_{k \in \mathcal{K}} \alpha_k\right) - \sum_{k \in \mathcal{K}} \alpha_k \theta_k Q \geq 0 \tag{3}$$

$$\alpha_k \theta_k Q - \mathcal{J}_{p_k}(\alpha_k) \geq 0 \quad \forall k \in \mathcal{K} \tag{4}$$

Then, we can get

$$\sum_{k \in \mathcal{K}} \alpha_k \theta_k \leq \frac{Q - \mathcal{J}_o(1 - \sum_{m \in \mathcal{K}} \alpha_m)}{Q} \tag{5}$$

$$\theta_k(\boldsymbol{\alpha}) \geq \frac{\mathcal{J}_{p_k}(\alpha_k)}{\alpha_k Q} \quad \forall k \in \mathcal{K} \tag{6}$$

Equations (5) and (6) define the feasible domain of the pair of $\boldsymbol{\alpha}$ and $\boldsymbol{\theta}$, which is denoted by $\boldsymbol{\Omega}$.

According to the utility functions defined in this section, we propose an incentive design framework in the following section.

## III. Incentive Design Framework

In this section, we propose an incentive design framework based on the utility functions, Eqs. (1) and (2), and the Nash bargaining solution. The bargaining problem studies how players share a jointly generated surplus. To formulate the bargaining problem, we need to define the disagreement point, denoted by $r$, which is a set of strategies that provide the lowest utility as expected by players if the bargaining breaks down. In our problem, the disagreement point is when commercial players refuse to participate in space infrastructure development and deployment, which means $r \equiv u(\boldsymbol{\alpha} = \mathbf{0}^{|\mathcal{K}| \times 1}, \boldsymbol{\theta}) = 0$. Setting the disagreement point to be zero also eliminates any potential negative incentives. We define $\mathcal{Z} = \{\text{coordinator}\} \cup \{\text{commercial players}\}$ as the set of all game players, indexed by $z$.

Nash proposed that a solution to the bargaining problem should satisfy four axioms [27]: 1) independence of equivalent utility transformations, 2) independence of irrelevant alternatives, 3) Pareto optimality, 4) symmetry. The unique solution that satisfies all these axioms is the particular utility vector that maximizes the Nash product $\prod_{z \in \mathcal{Z}}(u_z(\boldsymbol{\alpha}, \boldsymbol{\theta}) - r_z)$. This solution was originally proposed for cooperative bargaining [27], but has also been shown to be relevant in the non-cooperative bargaining game [28]. In our case, because the disagreement point $r_z = 0$ for all players $z$, the problem of interest becomes:

maximize

$$\prod_{z \in \mathcal{Z}} (u_z(\boldsymbol{\alpha}, \boldsymbol{\theta})) \tag{7}$$

subject to

$$(\boldsymbol{\alpha}, \boldsymbol{\theta}) \in \Omega \tag{8}$$

Note that this feasible region in Eq. (8) ensures the incentive is beneficial to both the coordinator and commercial players:

$$u_z(\boldsymbol{\alpha}, \boldsymbol{\theta}) \geq 0 \quad \forall z \in \mathcal{Z} \tag{9}$$

Depending on different mission planning circumstances, there are three different mission scenarios for decision variables to be discussed:

- *Scenario 1:* Optimal design: both $\boldsymbol{\alpha}$ and $\boldsymbol{\theta}$ are decision variables – This case corresponds to the mission planning scenario where the coordinator is trying to find the optimal design point without a pre-determined payload transportation assignment and the limitation of incentive budget.

- *Scenario 2:* Fixed demand allocation: $\boldsymbol{\alpha}$ is a given constant vector and $\boldsymbol{\theta}$ is the only decision variable – This case corresponds to the scenario where the coordinator has already determined the transportation task assignment for each player, the only problem is how to identify the optimal incentive.

- *Scenario 3:* Fixed budget allocation: $\boldsymbol{\alpha}$ is the only decision variable and $\boldsymbol{\theta}$ is a given constant vector – This case corresponds to the scenario where the coordinator has to decide the optimal space transportation task assignment given certain total mission demand under the limitation of a given incentive budget allocation.

In the following sections, we discuss the properties of each design scenario and approaches to solve the problem.

**A.** *Scenario 1*: **Optimal Design**

When both $\boldsymbol{\alpha}$ and $\boldsymbol{\theta}$ are decision variables, we can prove two important properties of the Nash bargaining solution ($\boldsymbol{\alpha}_{\text{NBS}}, \boldsymbol{\theta}_{\text{NBS}}$) with respect to the total utility and the utility distribution.

**Theorem 1.** *Given the utility definition in Eqs. (1) & (2) and zero disagreement points, for a multi-player incentive design problem with one coordinator, the Nash bargaining solution ($\boldsymbol{\alpha}_{\text{NBS}}, \boldsymbol{\theta}_{\text{NBS}}$) maximizes the total social welfare, which is the total utility of the coordinator and commercial players.*

$$\sum_{z \in \mathcal{Z}} u_z(\boldsymbol{\alpha}_{\text{NBS}}, \boldsymbol{\theta}_{\text{NBS}}) \geq \sum_{z \in \mathcal{Z}} u_z(\boldsymbol{\alpha}', \boldsymbol{\theta}') \quad \forall (\boldsymbol{\alpha}', \boldsymbol{\theta}') \in \Omega \tag{10}$$

*PROOF.* Based on the definition of utility functions, we know the total social welfare can be expressed as

$$\sum_{z \in Z} u_z(\pmb{\alpha}, \pmb{\theta}) = Q - J_o\left(1 - \sum_{k \in \mathcal{K}} \alpha_k\right) - \sum_{k \in \mathcal{K}} J_{p_k}(\alpha_k)$$

which is independent of the incentive coefficient, $\pmb{\theta}$. Namely, $\pmb{\theta}$ only influences the distribution of utilities but does not impact the total utility. We define $U(\pmb{\alpha}) = \sum_{z \in Z} u_z(\pmb{\alpha}, \pmb{\theta})$ as the total utility of all players (i.e., total social welfare).

Given an $\pmb{\alpha}$ value, the total utility of all players $U(\pmb{\alpha})$ is fixed, and the optimal utility distribution to maximize the Nash product under that constraint is achieved by making the utility of all players equal if that solution is feasible. Thus, we first show the feasibility of such a solution in the considered scenario; namely, for any specific positive participation coefficients $\pmb{\alpha}$, we can always find an incentive coefficient vector $\pmb{\theta}$ in the domain $\pmb{\Omega}$ to make the coordinator and commercial players achieve an equal utility, $u_o(\pmb{\alpha}, \pmb{\theta}) = u_{p_k}(\alpha_k, \theta_k) \; \forall k \in \mathcal{K}$. Then, we only need to prove that the Nash bargaining solution maximizes the social welfare (i.e., by varying $\pmb{\alpha}$) when the utility for each player is equal.

First, to make $u_o(\pmb{\alpha}, \pmb{\theta}) = u_{p_k}(\alpha_k, \theta_k) \; \forall k \in \mathcal{K}$, we have

$$u_o(\pmb{\alpha}, \pmb{\theta}) = u_{p_k}(\alpha_k, \theta_k) = Q - J_o\left(1 - \sum_{m \in \mathcal{K}} \alpha_m\right) - \sum_{m \in \mathcal{K}} \alpha_m \theta_m Q = \alpha_k \theta_k Q - J_{p_k}(\alpha_k) \quad \forall k \in \mathcal{K} \quad (11)$$

$$u_{p_k}(\alpha_k, \theta_k) = u_{p_m}(\alpha_m, \theta_m) = \alpha_k \theta_k Q - J_{p_k}(\alpha_k) = \alpha_m \theta_m Q - J_{p_m}(\alpha_m) \quad \forall m \neq k, \forall m, k \in \mathcal{K} \quad (12)$$

Combining Eqs. (11) and (12), we can get

$$\theta_k^*(\pmb{\alpha}) = \frac{Q - J_o(1 - \sum_{m \in \mathcal{K}} \alpha_m) + N J_{p_k}(\alpha_k) - \sum_{m \in \mathcal{K}} J_{p_m}(\alpha_m)}{N \alpha_k Q} \quad \forall k \in \mathcal{K}$$

Here, we show $(\pmb{\alpha}, \pmb{\theta}) \in \pmb{\Omega}$ using Eqs. (5) and (6). We know that $\theta_k^*(\pmb{\alpha})$ satisfies Eq. (5) because

$$\frac{Q - J_o(1 - \sum_{m \in \mathcal{K}} \alpha_m)}{Q} - \sum_{k \in \mathcal{K}} \alpha_k \theta_k^*(\pmb{\alpha}) = \frac{Q - J_o(1 - \sum_{m \in \mathcal{K}} \alpha_m) - \sum_{m \in \mathcal{K}} J_{p_m}(\alpha_m)}{NQ} = \frac{\sum_{z \in Z} u_z(\pmb{\alpha}, \pmb{\theta})}{NQ} \geq 0$$

We also know that $\theta_k^*(\pmb{\alpha})$ satisfies Eq. (6) because

$$\theta_k^*(\pmb{\alpha}) - \frac{J_{p_k}(\alpha_k)}{\alpha_k Q} = \frac{Q - J_o(1 - \sum_{m \in \mathcal{K}} \alpha_m) - \sum_{m \in \mathcal{K}} J_{p_m}(\alpha_m)}{N \alpha_k Q} = \frac{\sum_{z \in Z} u_z(\pmb{\alpha}, \pmb{\theta})}{N \alpha_k Q} \geq 0 \quad \forall k \in \mathcal{K}$$

Thus, for all $\pmb{\alpha}$, we have

$$(\pmb{\alpha}, \pmb{\theta}) \in \pmb{\Omega}$$

Next, we can show that the Nash bargaining solution gives the maximal social welfare when utilities are equal, which is equivalent to the statement that the Nash bargaining solution maximizes the total utility in general.

For the Nash bargaining solution, we have

$$\prod_{z\in Z} u_z(\boldsymbol{\alpha}_{\text{NBS}}, \boldsymbol{\theta}_{\text{NBS}}) \geq \prod_{z\in Z} u_z(\boldsymbol{\alpha}', \boldsymbol{\theta}') \quad \forall (\boldsymbol{\alpha}', \boldsymbol{\theta}') \in \Omega$$

When the utilities are equal, we have

$$u_o(\boldsymbol{\alpha}, \boldsymbol{\theta}) = u_{p_k}(\alpha_k, \theta_k) = \frac{Q - \mathcal{J}_o(1 - \sum_{m\in\mathcal{K}} \alpha_m) - \sum_{m\in\mathcal{K}} \mathcal{J}_{p_m}(\alpha_m)}{N} \quad \forall k \in \mathcal{K}$$

Thus, we get,

$$\frac{(Q - \mathcal{J}_o(1 - \sum_{m\in\mathcal{K}} \alpha_{m,\text{NBS}}) - \sum_{m\in\mathcal{K}} \mathcal{J}_{p_m}(\alpha_{m,\text{NBS}}))^N}{N^N} \geq \frac{(Q - \mathcal{J}_o(1 - \sum_{m\in\mathcal{K}} \alpha'_m) - \sum_{m\in\mathcal{K}} \mathcal{J}_{p_m}(\alpha'_m))^N}{N^N} \quad \forall \boldsymbol{\alpha}'$$

Because the total utility is defined to be nonnegative (i.e., $U(\boldsymbol{\alpha}) \geq 0$),

$$Q - \mathcal{J}_o\left(1 - \sum_{m\in\mathcal{K}} \alpha_{m,\text{NBS}}\right) - \sum_{m\in\mathcal{K}} \mathcal{J}_{p_m}(\alpha_{m,\text{NBS}}) \geq Q - \mathcal{J}_o\left(1 - \sum_{m\in\mathcal{K}} \alpha'_m\right) - \sum_{m\in\mathcal{K}} \mathcal{J}_{p_m}(\alpha'_m) \quad \forall \boldsymbol{\alpha}'$$

This inequality leads to

$$\sum_{z\in Z} u_z(\boldsymbol{\alpha}_{\text{NBS}}, \boldsymbol{\theta}_{\text{NBS}}) \geq \sum_{z\in Z} u_z(\boldsymbol{\alpha}', \boldsymbol{\theta}') \quad \forall (\boldsymbol{\alpha}', \boldsymbol{\theta}') \in \Omega$$

Therefore, we can claim that the Nash bargaining solution of the incentive design also maximizes social welfare.

∎

Besides the total utility obtained by the system, we also care about the fairness of the benefit distribution. The incentive can be designed focusing on fairness leveraging the maximin strategy. The purpose of the maximin strategy is to optimize the utility distribution among players. It achieves this goal by maximizing the minimum utility among players. For the Nash bargaining solution, we can also prove the following theorem regarding the utility distribution.

**Theorem 2.** *Given the utility definition in Eqs. (1) & (2) and zero disagreement points, for a multi-player incentive design problem with one coordinator, the Nash bargaining solution ($\boldsymbol{\alpha}_{\text{NBS}}, \boldsymbol{\theta}_{\text{NBS}}$) also maximizes the minimum utility among players.*

$$\min\{u_z(\boldsymbol{\alpha}_{\text{NBS}}, \boldsymbol{\theta}_{\text{NBS}}), z \in \mathcal{Z}\} \geq \min\{u_z(\boldsymbol{\alpha}', \boldsymbol{\theta}'), z \in \mathcal{Z}\} \quad \forall (\boldsymbol{\alpha}', \boldsymbol{\theta}') \in \Omega \tag{13}$$

*PROOF.* Because of the symmetry of the utility space in our problem, the symmetry axiom of the Nash bargaining solution shows

$$\min\{u_z(\boldsymbol{\alpha}_{\text{NBS}}, \boldsymbol{\theta}_{\text{NBS}}), z \in \mathcal{Z}\} = \frac{Q - \mathcal{J}_o(1 - \sum_{m\in\mathcal{K}} \alpha_m) - \sum_{m\in\mathcal{K}} \mathcal{J}_{p_m}(\alpha_m)}{N}$$

Thus, based on Theorem 1, we get

$$N \min\{u_z(\boldsymbol{\alpha}_{\text{NBS}}, \boldsymbol{\theta}_{\text{NBS}}), z \in \mathcal{Z}\} = \sum_{z\in Z} u_z(\boldsymbol{\alpha}_{\text{NBS}}, \boldsymbol{\theta}_{\text{NBS}}) \geq \sum_{z\in Z} u_z(\boldsymbol{\alpha}', \boldsymbol{\theta}') \geq N \min\{u_z(\boldsymbol{\alpha}', \boldsymbol{\theta}'), z \in \mathcal{Z}\}$$

∎

**Remark 1.** *Theorems 1-2 can be easily extended to the case with nonzero disagreement points by focusing on the surplus utility $u_z - r_z$ in place of the utility $u_z$ itself.*

Based on Theorem 1 proven above, we do not need to solve the original nonlinear formulation with the multiplicative objective function (i.e., Nash product) as shown in Eqs. (7)-(9). Instead, we propose an analytical approach and a numerical mixed-integer linear programming (MILP) approach as follows to solve the incentive design problem. In these two approaches, we have objective functions that are only linearly dependent on the mission cost functions $\mathcal{J}_o$ and $\mathcal{J}_p$.

*1. Analytical Approach*

If we know the analytical expression of mission cost functions $\mathcal{J}_o(\boldsymbol{\alpha})$ and $\mathcal{J}_p(\boldsymbol{\alpha})$, we can solve the incentive design problem analytically. Based on Theorem 1, we can identify the optimal $\boldsymbol{\alpha}_{\text{NBS}}$ for the Nash bargaining solution by maximizing the summation of total utilities,

$$\boldsymbol{\alpha}_{\text{NBS}} = \arg\max_{\boldsymbol{\alpha}} Q - \mathcal{J}_o\left(1 - \sum_{m \in \mathcal{K}} \alpha_m\right) - \sum_{m \in \mathcal{K}} \mathcal{J}_{p_m}(\alpha_m) \qquad (14)$$

Then, the optimal $\boldsymbol{\theta}_{NBS}$ can be determined by letting utilities be equal, $u_o(\boldsymbol{\alpha}, \boldsymbol{\theta}) = u_{p_k}(\alpha_k, \theta_k) \ \forall k \in \mathcal{K}$, which gives

$$\theta_{k,\text{NBS}} = \frac{Q - \mathcal{J}_o(1 - \sum_{m \in \mathcal{K}} \alpha_{m,\text{NBS}}) + N\mathcal{J}_{p_k}(\alpha_{k,\text{NBS}}) - \sum_{m \in \mathcal{K}} \mathcal{J}_{p_m}(\alpha_{m,\text{NBS}})}{N\alpha_{k,\text{NBS}}Q} \qquad \forall k \in \mathcal{K} \qquad (15)$$

The analytical solution is convenient if we know expressions of mission cost functions explicitly. Some studies have proposed approximate mission cost functions [25, 26]. However, the application of these cost functions is limited to certain mission scenarios to achieve accurate approximation. We mainly rely on the numerical approach to solve the incentive design problem for most general mission planning problems, which will be discussed in the next subsection.

*2. Numerical Approach*

If we do not know the accurate analytical expression of mission cost functions $\mathcal{J}_o(\boldsymbol{\alpha})$ and $\mathcal{J}_p(\boldsymbol{\alpha})$, we can solve the incentive design problem numerically as a MILP problem leveraging the network-based space logistics optimization model [8]. In this problem, the nonlinear objective function in Eq. (7) is converted into a linear one using Theorem 1.

The proposed formulation considers space logistics as a multi-commodity network flow problem. In the network, nodes represent orbits or planets; arcs represent space flight trajectories; crew, propellant, instruments, spacecraft, and all other payloads are considered as commodities flowing along arcs, as shown in Fig. 2. A time dimension is introduced to take into account time steps for dynamic mission planning. There are mainly two types of arcs: 1) transportation arcs that connect different nodes at different time steps to represent space flights; 2) holdover arcs that connect the same node at different time steps to represent in-orbit or surface space mission operations after space infrastructure deployment.

Consider a network graph defined by a set of arcs $\mathcal{A} = \{\mathcal{Z}, \mathcal{V}, \mathcal{N}, \mathcal{T}\}$, which includes a player set $\mathcal{Z}$ (index: $z$), a spacecraft set $\mathcal{V}$ (index: $v$), a node set $\mathcal{N}$ (index: $i$ and $j$), and a time step index set $\mathcal{T}$ (index: $t$). We define a commodity flow variable $\boldsymbol{x}_{zvijt}$ as the decision variable for space mission planning, denoting the commodity flow from node $i$ to node $j$ at time step $t$ using spacecraft $v$ performed by player $z$. If we denote the set of commodities as $\mathcal{C}$, this $\boldsymbol{x}_{zvijt}$ is a $|\mathcal{C}| \times 1$ vector, and its elements include both continuous and discrete variables depending on the commodity type. For example, a continuous variable is used for the mass of payload or propellant; whereas a discrete variable is used for the number of crew members or spacecraft. We define $\mathcal{C}_\text{C}$ and $\mathcal{C}_\text{D}$ to denote the sets of commodities that are considered as continuous and discrete variables, respectively. Then, the commodity flow variable can be expressed as

$$\boldsymbol{x}_{zvijt} = \begin{bmatrix} \boldsymbol{x}_\text{C} \\ \boldsymbol{x}_\text{D} \end{bmatrix}_{zvijt}$$

where $\boldsymbol{x}_\text{C}$ is a $|\mathcal{C}_\text{C}| \times 1$ vector for continuous commodities and $\boldsymbol{x}_\text{D}$ is a $|\mathcal{C}_\text{D}| \times 1$ vector for discrete commodities.

To measure the space mission planning performance, we also define a cost coefficient $\boldsymbol{c}_{zvijt}$. The mission demands or supplies are determined through a demand vector $\boldsymbol{d}_{zit}$, denoting the mission demand or supply for player $z$ in node $i$ at time $t$. Mission demands are negative, and mission supplies are positive. As there are $|\mathcal{C}|$ types of commodities, $\boldsymbol{c}_{zvijt}$ and $\boldsymbol{d}_{zit}$ are both $|\mathcal{C}| \times 1$ vectors.

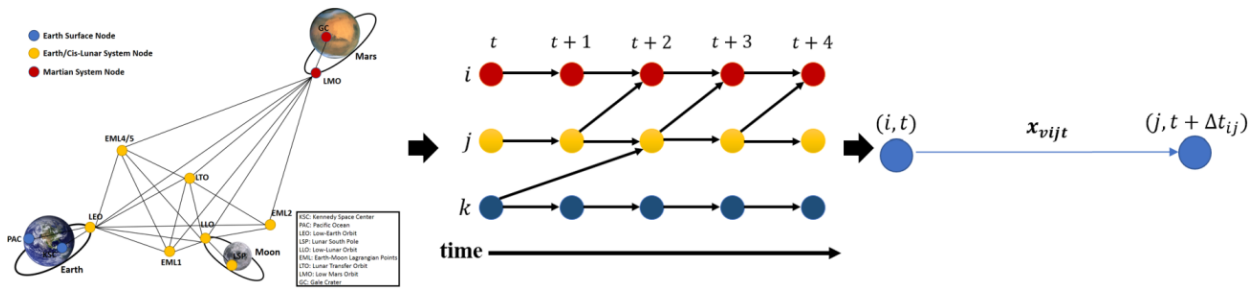

**Fig. 2 Space logistics network model [6, 11]**

Besides the parameters defined above, we also need to define the following parameters for the formulation:

$\Delta t_{ij}$ = time of flight.

$G_{zvij}$ = commodity transformation matrix.

$H_{zvij}$ = concurrency constraint matrix.

$W_{zij}$ = mission time windows.

Based on the aforementioned notations, the problem formulation for social welfare maximization is shown as follows:

maximize

$$U(\boldsymbol{\alpha}) = \sum_{z \in Z} u_z(\boldsymbol{\alpha}_{\text{NBS}}, \boldsymbol{\theta}_{\text{NBS}}) = Q - \sum_{(z,v,i,j,t) \in \mathcal{A}} \boldsymbol{c}_{zvijt}^T \boldsymbol{x}_{zvijt} \tag{16}$$

subject to

$$\sum_{(z,v,i,j,t) \in \mathcal{A}} \boldsymbol{c}_{zvijt}^T \boldsymbol{x}_{zvijt} \leq Q \tag{17}$$

$$\sum_{(v,j):(z,v,i,j,t) \in \mathcal{A}} \boldsymbol{x}_{zvijt} - \sum_{(v,j):(z,v,j,i,t) \in \mathcal{A}} G_{zvji} \boldsymbol{x}_{zvji(t-\Delta t_{ji})} \leq \boldsymbol{d}_{zit} + \boldsymbol{d}'_{zit}(\boldsymbol{\alpha}) \quad \forall z \in Z \;\; \forall i \in \mathcal{N} \;\; \forall t \in \mathcal{T} \tag{18}$$

$$H_{zvij} \boldsymbol{x}_{zvijt} \leq \boldsymbol{0}_{l \times 1} \quad \forall (z,v,i,j,t) \in \mathcal{A} \tag{19}$$

$$\begin{cases} \boldsymbol{x}_{zvijt} \geq \boldsymbol{0}_{|\mathcal{C}| \times 1} & \text{if } t \in W_{zij} \\ \boldsymbol{x}_{zvijt} = \boldsymbol{0}_{|\mathcal{C}| \times 1} & \text{otherwise} \end{cases} \quad \forall (z,v,i,j,t) \in \mathcal{A} \tag{20}$$

$$\boldsymbol{x}_{zvijt} = \begin{bmatrix} \boldsymbol{x}_{\text{C}} \\ \boldsymbol{x}_{\text{D}} \end{bmatrix}_{zvijt}, \boldsymbol{x}_{\text{C}} \in \mathbb{R}_{\geq 0}^{|\mathcal{C}_{\text{C}}| \times 1}, \boldsymbol{x}_{\text{D}} \in \mathbb{Z}_{\geq 0}^{|\mathcal{C}_{\text{D}}| \times 1} \quad \forall (z,v,i,j,t) \in \mathcal{A} \tag{21}$$

In this formulation, Eq. (16) is the objective function that maximizes social welfare. After obtaining the optimal value of $\boldsymbol{\alpha}$ from this optimization, we can calculate the optimal $\boldsymbol{\theta}$ using the Eq. (15) introduced earlier. These optimal $\boldsymbol{\alpha}$ and $\boldsymbol{\theta}$ lead to the Nash bargaining solution by Theorem 1. The term $\sum_{(z,v,i,j,t) \in \mathcal{A}} \boldsymbol{c}_{zvijt}^T \boldsymbol{x}_{zvijt}$ is the total mission cost to complete the infrastructure deployment mission, which is equivalent to $\mathcal{J}_o(1 - \sum_{m \in \mathcal{K}} \alpha_m) + \sum_{m \in \mathcal{K}} \mathcal{J}_{p_m}(\alpha_m)$. Equation (17) makes sure that the total mission cost is always smaller or equal to the baseline mission cost, which creates a non-negative total utility for the system. Equation (18) is the commodity mass balance constraint. It guarantees that the commodity inflow is always larger or equal to the commodity outflow plus the mission demands. The second term $G_{zvji} \boldsymbol{x}_{zvji(t-\Delta t_{ji})}$ in this constraint represents commodity transformations during space flights or operations, including propellant burning, crew consumptions, and ISRU resource productions. For the detailed settings of the transformation matrix $G_{zvij}$, please refer to Ref. [8]. In this constraint, the demand vector $\boldsymbol{d}_{zit}$ denotes the original mission demand depending on each player's mission goal. Another demand vector $\boldsymbol{d}'_{zit}(\boldsymbol{\alpha})$ denotes additional infrastructure deployment demand determined based on the participation coefficient $\alpha$. Equation (19) is the concurrency constraint that defines the commodity flow upper bound typically determined by the spacecraft payload and propellant capacities. The index $l$ is for the type of concurrency constraints considered. Equation (20) is the time bound. It defines the space flight and operation time windows. Only when the time windows are open, is the space transportation along the arc permitted.

This linear formulation can solve the incentive design optimization problem in a computationally efficient way without relying on the original nonlinear optimization formulation in Eqs. (7)-(9).

B. *Scenario 2*: **Fixed Demand Allocation**

If $\alpha$ is a given constant vector and $\theta$ is the only decision variable, the coordinator has already determined the deployment task assignment for each player, the only problem is how to optimize the incentive.

This scenario is a special case of Scenario 1. Theorem 1 in Scenario 1 is still valid because $\alpha$ only determines the total utility of the system, while $\theta$ determines the utility distribution among players. Thus, for any given $\alpha^*$, we can always determine the optimal $\theta_{\text{NBS}}$ by the same approach as Scenario 1,

$$\theta_{k,\text{NBS}} = \frac{Q - \mathcal{J}_o(1 - \sum_{m \in \mathcal{K}} \alpha_m^*) + N\mathcal{J}_{p_k}(\alpha_k^*) - \sum_{m \in \mathcal{K}} \mathcal{J}_{p_m}(\alpha_m^*)}{N\alpha_k^* Q} \quad \forall k \in \mathcal{K} \tag{22}$$

Similarly, players' utilities and mission costs can be solved either through the analytical approach or the numerical formulation introduced in Sec. III.A.

C. *Scenario 3*: **Fixed Budget Allocation**

If $\alpha$ is the only decision variable and $\theta$ is a given constant vector, the coordinator has to make the deployment task assignment under the incentive budget limitation. Assuming that the given incentive coefficient is $\theta^*$, we can define the incentive budget for the commercial player $k$ as $R_k = \theta_k^* Q$. Then, the utilities for the coordinator and the commercial player $k$ can be rewritten as,

$$u_o(\alpha, \theta) = Q - \mathcal{J}_o\left(1 - \sum_{k \in \mathcal{K}} \alpha_k\right) - \sum_{k \in \mathcal{K}} \alpha_k R_k \tag{23}$$

$$u_{p_k}(\alpha_k, \theta_k) = \alpha_k R_k - \mathcal{J}_{p_k}(\alpha_k) \tag{24}$$

Now that the participation coefficient vector $\alpha$ is the only decision variable, we can see that Theorem 1 no longer holds. This is because the solution we used to derive these theorems (i.e., equating the utilities of all players) is not feasible as $\boldsymbol{\theta}$ cannot be freely chosen. A fixed incentive coefficient vector $\boldsymbol{\theta}$ cannot guarantee equal utilities for players. We have to solve the original nonlinear optimization formulation, Eqs. (7)-(9), to identify the optimal $\boldsymbol{\alpha}^*$,

$$\boldsymbol{\alpha}^* = \arg\max_{\boldsymbol{\alpha}} \prod_{z \in Z} u_z(\boldsymbol{\alpha}) \tag{25}$$

By implementing the network-based space logistics optimization model, we can express the utility $u_z(\boldsymbol{\alpha})$ as a linear function of $\boldsymbol{\alpha}$. The original nonlinear problem is then converted into a mixed-integer nonlinear programming formulation as follows:

maximize

$$\prod_{z \in Z} u_z(\boldsymbol{\alpha}) \tag{26}$$

subject to

Eqs. (17)-(21)

## IV. Lunar Exploration Case Study

In this section, we perform case studies of a lunar exploration campaign to demonstrate the performance of the proposed incentive design framework and the value of commercialization in space infrastructure development and deployment. First, in Sec. IV.A, we introduce the problem settings and mission scenarios for the analysis. Then, in Sec. IV.B, we perform numerical experiments to discuss the impact of the proposed incentive design framework. We also conduct a sensitivity analysis to evaluate the impact of deployment demands and ISRU productivities.

### A. Problem Settings

The coordinator (player 0) has annual infrastructure deployment demands; whereas commercial players (player 1,2,3,…) focus on their own lunar exploration mission. For simplicity, we assume the commercial players are the only players who have the ability to develop and deploy ISRU systems on the lunar surface to support space transportation. Commercial players can participate in the infrastructure deployment mission if the coordinator provides enough incentives.

We consider a cislunar transportation system, containing Earth, low-Earth orbit (LEO), Earth-Moon Lagrange point 1 (EML1), and the Moon. It is a four-node transportation network model, as shown in Fig. 3, where the discrete time-expanded network model uses one day as one time step. The space transportation $\Delta V$ is also shown in Fig. 3.

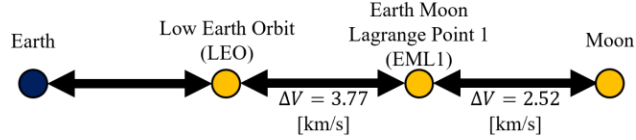

**Fig. 3 Cislunar Transportation Network Model**

For simplicity, we assume all players use identical spacecraft to conduct the space infrastructure deployment mission. The spacecraft has RL-10 rocket engines using liquid hydrogen and liquid oxygen (LH2/LOX) as the propellant. Based on the spacecraft sizing of the Advanced Cryogenic Evolved Stage (ACES) [29] spacecraft and the Centaur [30] spacecraft, we assume the transportation spacecraft considered has a dry mass of 6,000 kg and an inert mass fraction of 0.1. We assume the player may also have deployed ISRU systems in advance to support the infrastructure deployment in the transportation system. The ISRU system is assumed as a water electrolysis ISRU, which generates oxygen and hydrogen from water. The baseline productivity of the ISRU is 5 kg water/yr/kg plant [11], meaning that a 1 kg ISRU plant can electrolyze 5 kg water per year on the lunar surface. We also take into account the ISRU system maintenance, which requires maintenance spare resupply annually. The mass of maintenance spares needed is equivalent to 5% of the ISRU system mass [11]. The space mission operation parameters and assumptions are summarized in Table 1.

**Table 1 Mission parameters and assumptions**

| Parameters | Assumed value |
|---|---|
| S/c propellant capacity, kg | 54,000 [29, 30] |
| S/c structure mass, kg | 6,000 [29, 30] |
| Propellant type | LH2/LOX |
| Propellant $I_{sp}$, s | 420 |
| Propellant $O_2$:$H_2$ ratio | 5.5 |
| Water ISRU productivity, kg $H_2O$/ yr/ kg plant | 5 [11] |
| ISRU maintenance, /yr | 5% plant mass [11] |

To analyze the space mission cost, we also need a space mission cost model. This paper uses the same mission cost model proposed by Chen et al. in a multi-actor space commercialization study [18]. It includes rocket launch cost to LEO, spacecraft manufacturing cost, space flight operation cost, and LH2/LOX propellant price on Earth. The cislunar transportation cost model is listed in Table 2.

**Table 2 Cislunar transportation cost model [18]**

| Parameter | Assumed value |
|---|---|
| Rocket Launch cost, /kg | $3,500 |
| Spacecraft manufacturing, /spacecraft | $148M |
| Spacecraft operation, /flight | $1M |
| LH2 price on Earth, /kg | $5.94 |
| LO2 price on Earth, /kg | $0.09 |

At the beginning of the space mission, we assume that each commercial player has already built a lunar base with a 10 metric ton (MT) water ISRU system. For the nominal mission scenario, the coordinator plans to deploy 30 MT lunar habitats and infrastructures to the Moon every year. The nominal mission demands and supplies for this case study are shown in Table 3. The baseline mission cost is the value that the coordinator completes the mission independently.

**Table 3 Mission demands and supplies**

| Actors | Payload Type | Node | Time, day | Demand, MT |
|---|---|---|---|---|
| Coordinator | Infrastructure | Earth | 360 (repeat annually) | +30 |
|  | Infrastructure | Moon | 360 (repeat annually) | -30 |
|  | Propellant | Earth | All the time | +∞ |
| Commercial Player | ISRU plant, maintenance, and propellant | Earth | All the time | +∞ |
|  | ISRU plant | Moon | 0 | +10 |

**B. Numerical Experiments**

In this section, we first conduct a detailed analysis of a two-player case, with one coordinator and one player. We perform a sensitivity analysis to examine the relationships among mission demands by the coordinator, the size of the ISRU plant deployed by the commercial player, and the design spaces of decision variables, $\alpha$ and $\theta$. In addition, we extend the analysis to a three-player case, where we analyze the performance of the incentive framework in generating effective development strategies for one coordinator and two commercial players.

*1. Two-Player Case: Nominal Case*

For the nominal mission scenario, we assume that the coordinator (i.e., player 0) plans to complete two consecutive lunar transportation missions. Except for the existing ISRU plant owned by the commercial player, no players are allowed to deploy extra ISRU plants during the mission. We first calculate the baseline mission cost for this problem through the network-based space logistics model, which gives us $2,058M.

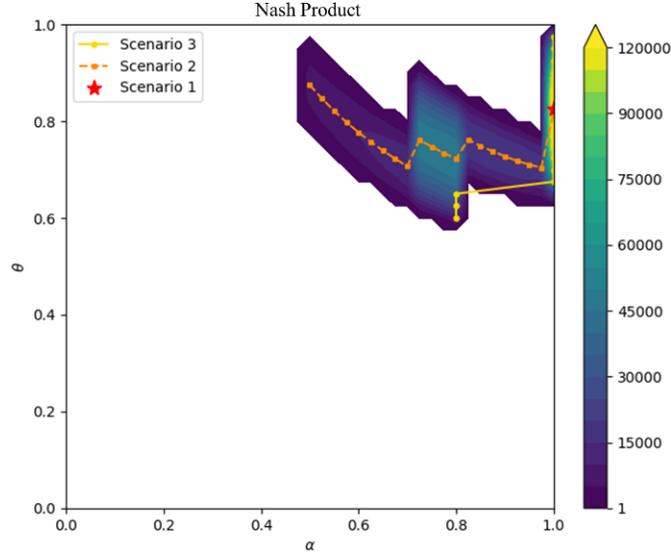

**Fig. 4 Nash Product Contour Plot and Optimal Solutions.**

The incentive design space contour and the optimal solutions for Scenario 1, 2, and 3 are shown in Fig. 4. Scenario 1 is the best solution within this contour plot; Scenario 2 is the best solution within this same plot given a value of $\alpha$; and Scenario 3 is the best solution within this same plot given the value of $\theta$. In this figure, the upper right corner means that the commercial player completes most of the deployment task, and the coordinator pays a relatively high incentive comparing with its own baseline mission cost. In the design space, moving left means less deployment task is assigned to the commercial player, while moving down means the coordinator pays lower incentive per unit mass. For this nominal mission scenario, the global optimal design point appears when the participation coefficient $\alpha$ is equal to 1, as shown as the red star in the figure. As expected, this Scenario 1 solution corresponds to the best one among the feasible solutions of Scenario 2 for all possible $\alpha$, which is also the best one among the feasible solutions of Scenario 3 for all possible $\theta$.

We first look into the impact of $\alpha$. From Fig. 4, we can find that the feasible incentive design space appears when $\alpha \geq \sim 0.5$. The feasible design region means that by leveraging the incentive design mechanism, the coordinator can complete its space infrastructure deployment mission at a lower total expense, while at the same time, the commercial player also receives profit through the difference between the incentive from the coordinator and its mission cost based on its existing infrastructure. (Note that it does not mean the mission is always feasible for any incentive coefficient $\theta$. It means that we can always find feasible and optimal $\theta$ when we are given an $\alpha$ in the domain.)

To further look into the feasible region of $\alpha$, we analyze the performance variation with $\alpha$ using the formulation in Scenario 2. The total expense comparison for this nominal mission between the baseline (i.e., in which the coordinator completes the mission independently) and the case with an incentive design implemented is shown in Fig. 5. Also, the optimal incentive paid to the commercial player and the total utility with respect to $\alpha$ is shown in Fig. 6.

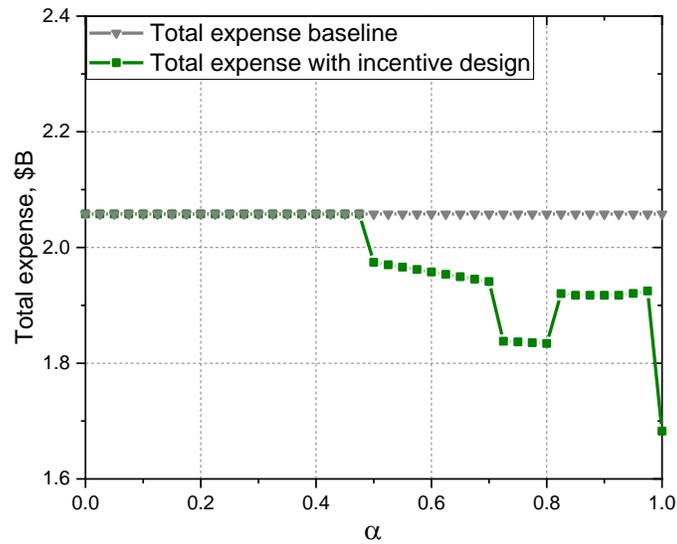

**Fig. 5 Total Expense Comparison (Scenario 2)**

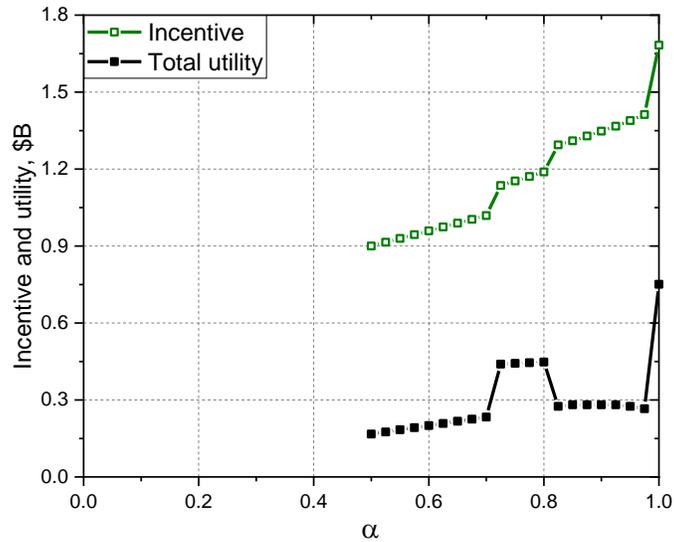

**Fig. 6 Total Utility and Incentive (Scenario 2)**

From Fig. 5, we can find that within the feasible design space of $\alpha$, the total expense of the coordinator after implementing the incentive design is lower than the baseline. Note that for the baseline case, the total expense is equal to the coordinator's mission cost; whereas, for the incentive design case, the total coordinator's expense is the summation of the coordinator's mission cost and the incentive paid to the commercial player. This result shows the value of space commercialization. By leveraging supporting infrastructures operated by other commercial players, the coordinator can reduce its own total expense to complete the deployment mission, even considering the additional incentive cost. The commercial player receives profit at the same time. Thus, the resulting solution is mutually beneficial for the government and the industry.

To better interpret the utility variation in the proposed model, for any given $\alpha$, we can also examine the mission costs of players 0 and 1 separately; note that these costs are calculated by the validated space mission design methods in the literature [6-8]. In this study, the surplus utilities come from the mission cost savings because the coordinator takes advantage of the power of commercial players. Therefore, based on the mission costs generated using traditional mission design methods, we are able to interpret the behavior and variation of the utility and the design space.

The mission cost of each player and the total transportation cost for any given $\alpha$ are shown in Fig. 7. Results in Fig. 7 explain why the feasible incentive design space appears when $\alpha \geq \sim 0.5$ for the considered mission scenario. This is the region where the total cost of space transportation starts to be smaller than the baseline cost (i.e., the total cost when $\alpha = 0$); this suggests that the two players can achieve more savings by working together than working independently (i.e., the utility is positive).

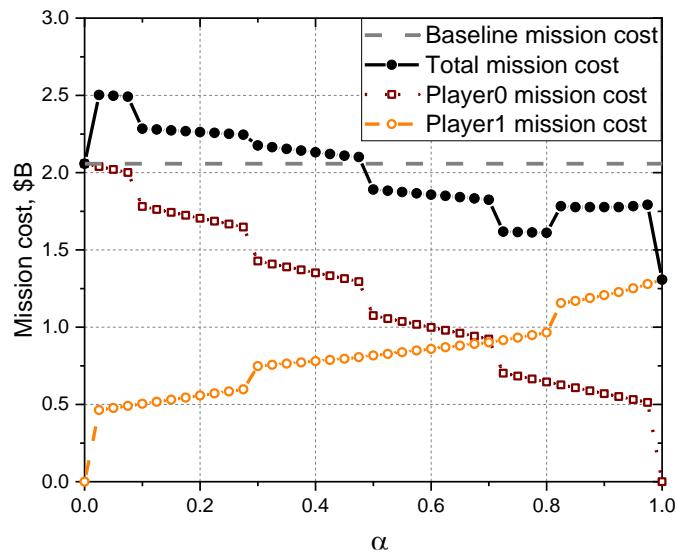

**Fig. 7 Mission Costs**

Moreover, it is also interesting to observe how the optimal value of $\theta$ changes as $\alpha$ changes. As can be observed in Scenario 2 in Fig. 4, as the increase of $\alpha$, the optimal $\theta$ for the Nash bargaining solution is decreasing. However, there are three sudden changes in optimal $\theta$ values during this process. These points correspond to three mission cost sudden drops or increases for players, as shown in Fig. 7. These dramatic changes in the transportation mission cost are initially caused by spacecraft because we consider the number of spacecraft as an integer variable, and the spacecraft manufacturing cost is also a significant portion of the transportation cost. When one more or less spacecraft is needed for transportation, there will be a sudden change in the mission cost. These changes are also reflected in the total expense of the mission in Fig. 5 and the total utility and incentive of the system in Fig. 6. The direction of the sudden change comes from the definition of utilities. We know that the incentive required is positively correlated to the cost of the commercial player (i.e., player 1) but negatively correlated to the cost of the coordinator (i.e., player 0). This explains that when there is a sudden decrease in the cost of the coordinator or a sudden increase in the cost of the commercial player, there will be a sudden increase in the value of optimal $\theta$. The choice of the optimal $\theta$ is the balance between the coordinator's and the commercial player's mission costs.

*2. Two-Player Case: Sensitivity Analysis*

We conduct a sensitivity analysis to analyze how the optimal solutions change with different parameter values. We consider the impact of the payload deployment demand and the size of the ISRU plant deployed in advance on the same nominal two-player mission scenario.

First, we consider a fixed 10 MT ISRU plant deployed in advance and only vary the space mission deployment demand. The incentive design space contours are shown in Fig. 8.

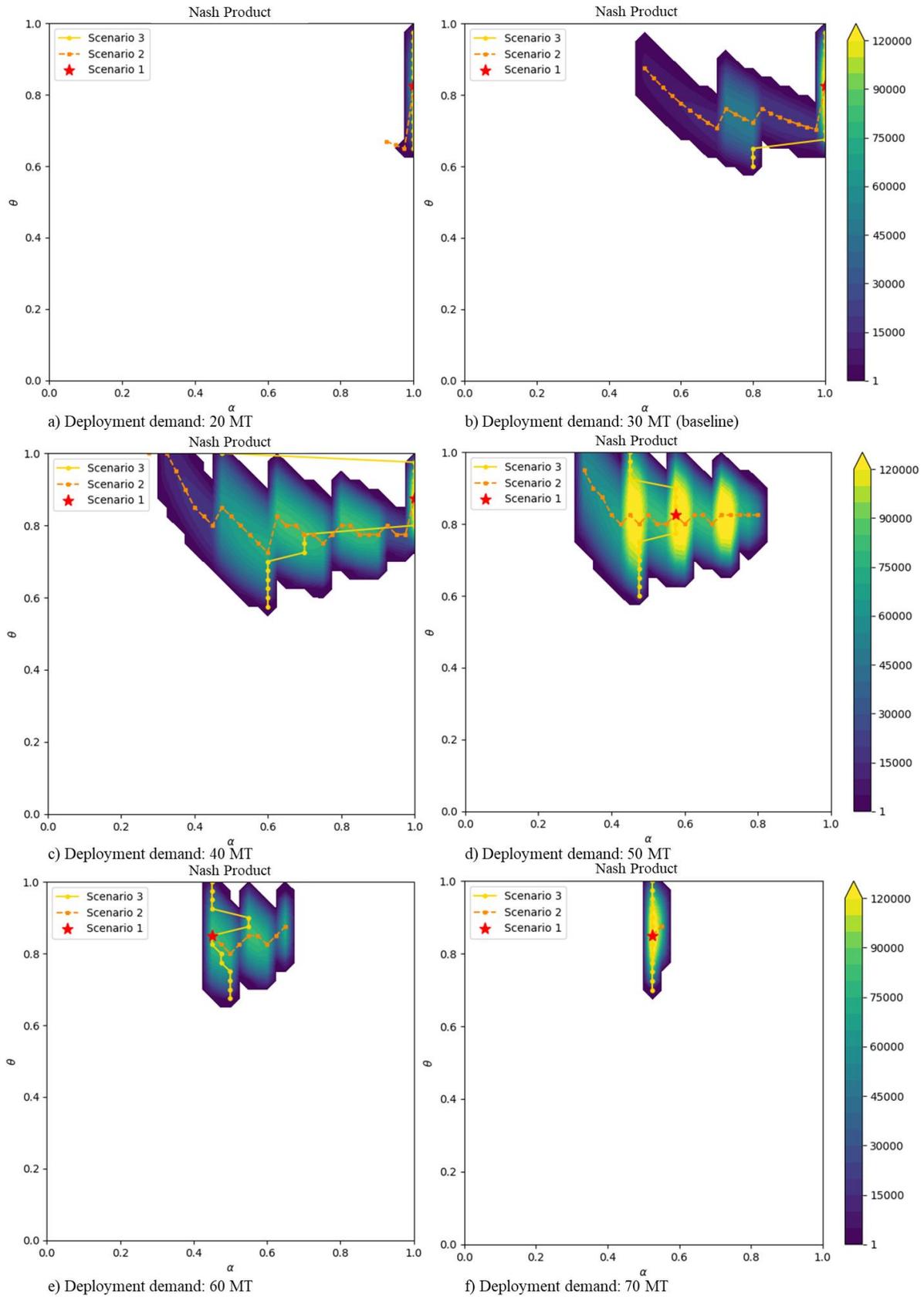

**Fig. 8 Incentive Design Space with Different Deployment Demands**

We can find that the mission demand determines the feasible domain of the participation coefficient $\alpha$ for the incentive design. When the deployment demand is small (i.e., 20 MT), the incentive design is only feasible when $\alpha$ is close to 1. As the demand increases, the feasible range for $\alpha$ expands from 1 to lower values. However, when the deployment demand is so large that no single player can complete the mission independently, any value that is close to 0 or 1 becomes infeasible. Eventually, the design space for $\alpha$ will converge to a constant value, which depends on the maximum transportation capacity of players. In this case, we assume that both players have two identical spacecraft at the beginning of each mission that provides a similar maximum transportation capacity. The maximum transportation capacities of the two players are not exactly equal because the ISRU plant operated by the commercial player can slightly increase its transportation ability though the impact is not significant. Thus, the optimal $\alpha$ converges to a value close to 0.5, as shown in Fig. 8, f).

The optimal $\alpha^*$ (i.e., the optimal solution of Scenario 1) also varies similarly to the feasible design space. For our case study, because of the existence of the ISRU plant, the commercial player can always complete the transportation mission at a lower cost than the coordinator. Therefore, when the deployment demand is low, such that a single player can complete the entire mission, the best strategy for the coordinator is to fully rely on the commercial player to complete the entire mission, which makes $\alpha^*$ to be 1. However, when the deployment demand exceeds a single player's capacity, $\alpha^*$ starts to converge to a value based on each player's maximum transportation capacity. Both the commercial player and the coordinator need to take over some portion of the deployment demand.

Next, we can assume the deployment demand is fixed to the nominal case of 30 MT, and only vary the size of ISRU plant deployment in advance. We can evaluate the change of design space with respect to the ISRU plant. The result is shown in Fig. 9.

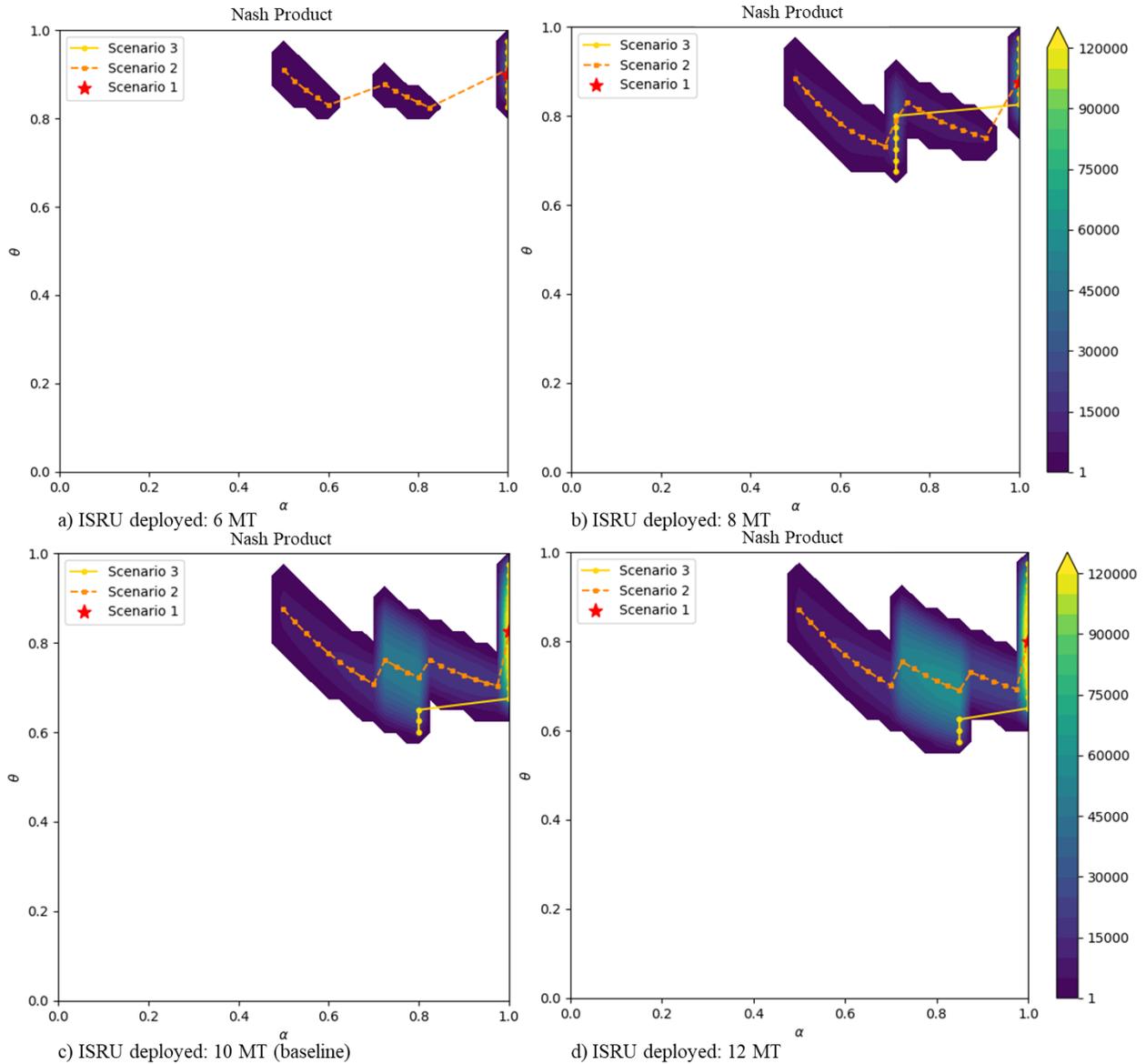

**Fig. 9 Incentive Design Space with Different Size of ISRU System Deployed**

It is shown that the size of the ISRU plant has a limited impact on the feasible domain of $\alpha$. Some portions of $\alpha$ may become infeasible in the domain $\sim 0.5 \leq \alpha \leq 1$, especially when the size of the ISRU is small. These infeasible gaps are caused by the sudden increase of total mission cost from spacecraft that eliminates positive system utility. But this impact becomes negligible when the ISRU becomes more productive (i.e., larger size).

However, the size of the ISRU plant directly determines the width of the feasible domain for the incentive coefficient $\theta$. When a larger ISRU plant is deployed, the commercial player can achieve considerable mission cost savings, which leads to a larger total utility. This result shows that a larger ISRU plant enables the decision-maker a more extensive design space for the incentive.

*3. Three-Player Case*

The proposed framework can also conduct performance analyses on missions with multiple commercial players. In this section, we consider a mission scenario with one coordinator and two commercial players. The mission planning assumptions and parameters are the same as the nominal case.

First, we consider a case where two commercial players have identical spacecraft and the same amount of ISRU plant deployed in advance. By changing the total deployment demand, we can observe the variation of the design space in the task assignment between two commercial players. The result is shown in Fig. 10. Note that this contour plot only shows $\alpha_1$ vs. $\alpha_2$, where $\boldsymbol{\theta}$ is always optimized to maximize the Nash products using the Eq. (15). The optimal points in the figure (labelled as red stars) are global optimal, which are also the solutions of Scenario 1. Because two commercial players are identical in this case, the contour plot should also be symmetric, and there may be two symmetric optimal points. The entire contour plot shows solutions to Scenario 2. Scenario 3 is not sghown as $\boldsymbol{\theta}$ is optimized. The dash lines indicate the optimal $\alpha$ for any given $\alpha$ of the other commercial player, which means the optimal $\alpha_1^*$ for any given $\alpha_2$ and the optimal $\alpha_2^*$ for any given $\alpha_1$.

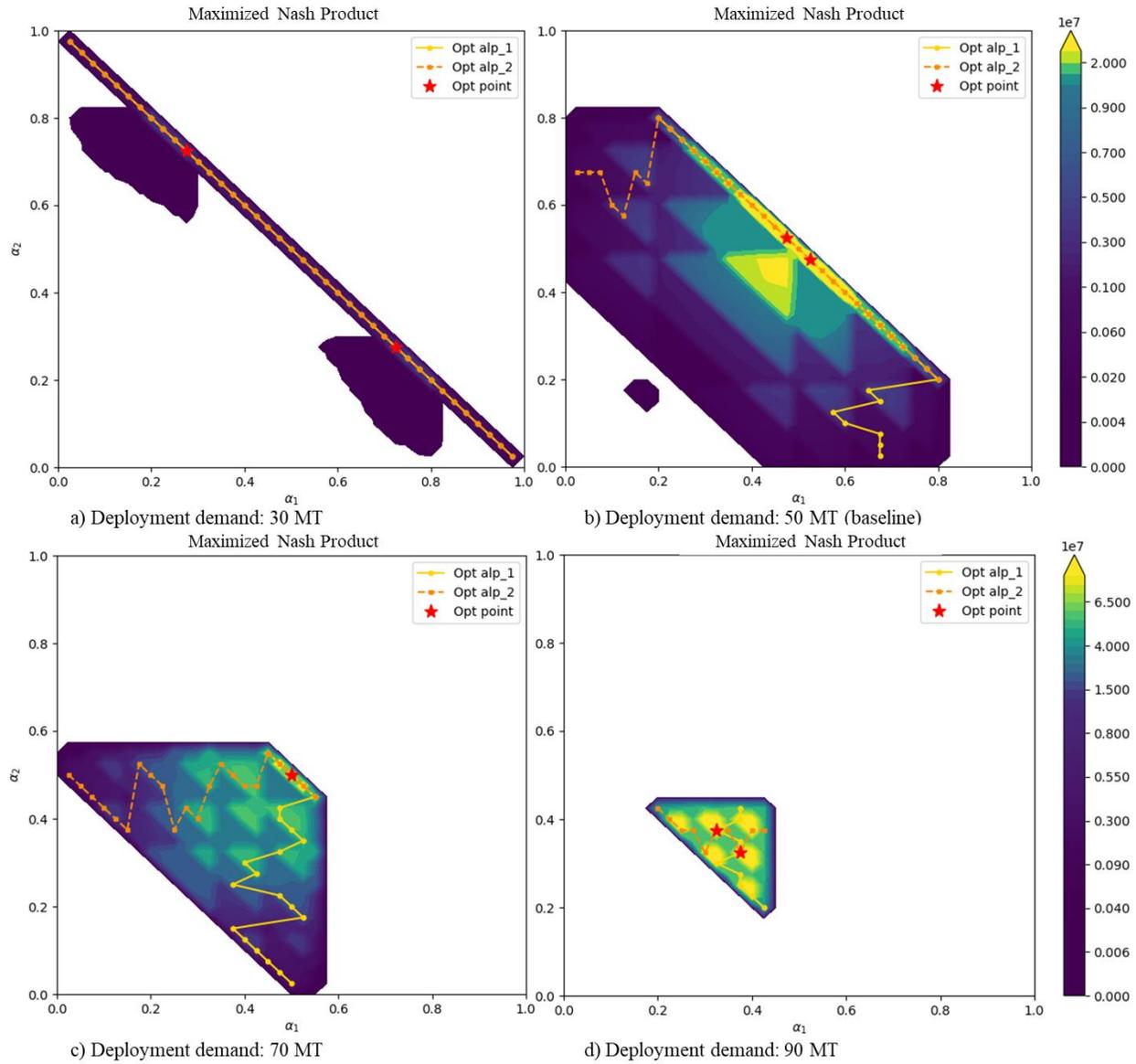

**Fig. 10 Three-Player Design Space with Different Deployment Demands**

The same trend can be observed as that for $\alpha$ value in Fig. 8. Similar to the two-player mission scenario, for this three-player mission, when the deployment demand is low, the optimal points only appear when $\alpha_1 + \alpha_2 \approx 1$. As the deployment demand increases, the total payload starts to exceed the transportation capability of commercial players. The coordinator starts to complement some transportation capacity.

Next, we consider a case where the capabilities of each commercial player is not equal to each other. Suppose we fix the total deployment demand to 50MT and vary the pre-deployed ISRU plant and the number of spacecraft operated by player 2. In that case, we can evaluate the imbalance in task assignment design space caused by the different sizes of the ISRU system and the transportation capacity. The result is shown in Fig. 11. For commercial player 1, the pre-deployed ISRU plant is fixed to 10 MT, and the number of spacecraft is fixed to 2 per mission as in the nominal mission scenario.

From the results, we can find that enlarging the ISRU plant and launching more spacecraft bring player 2 additional transportation capacity. In both cases, player 2 becomes capable to complete more deployment tasks. However, the number of spacecraft may constrain the potential improvement achieved by enlarging the ISRU plant. When player 2 only has one spacecraft, further increasing the size of ISRU cannot provide many benefits in the mission. This example demonstrates that the proposed framework is not only helpful for the coordinator to conduct incentive design but also beneficial to commercial players to identify effective development strategies in the commercialization environment.

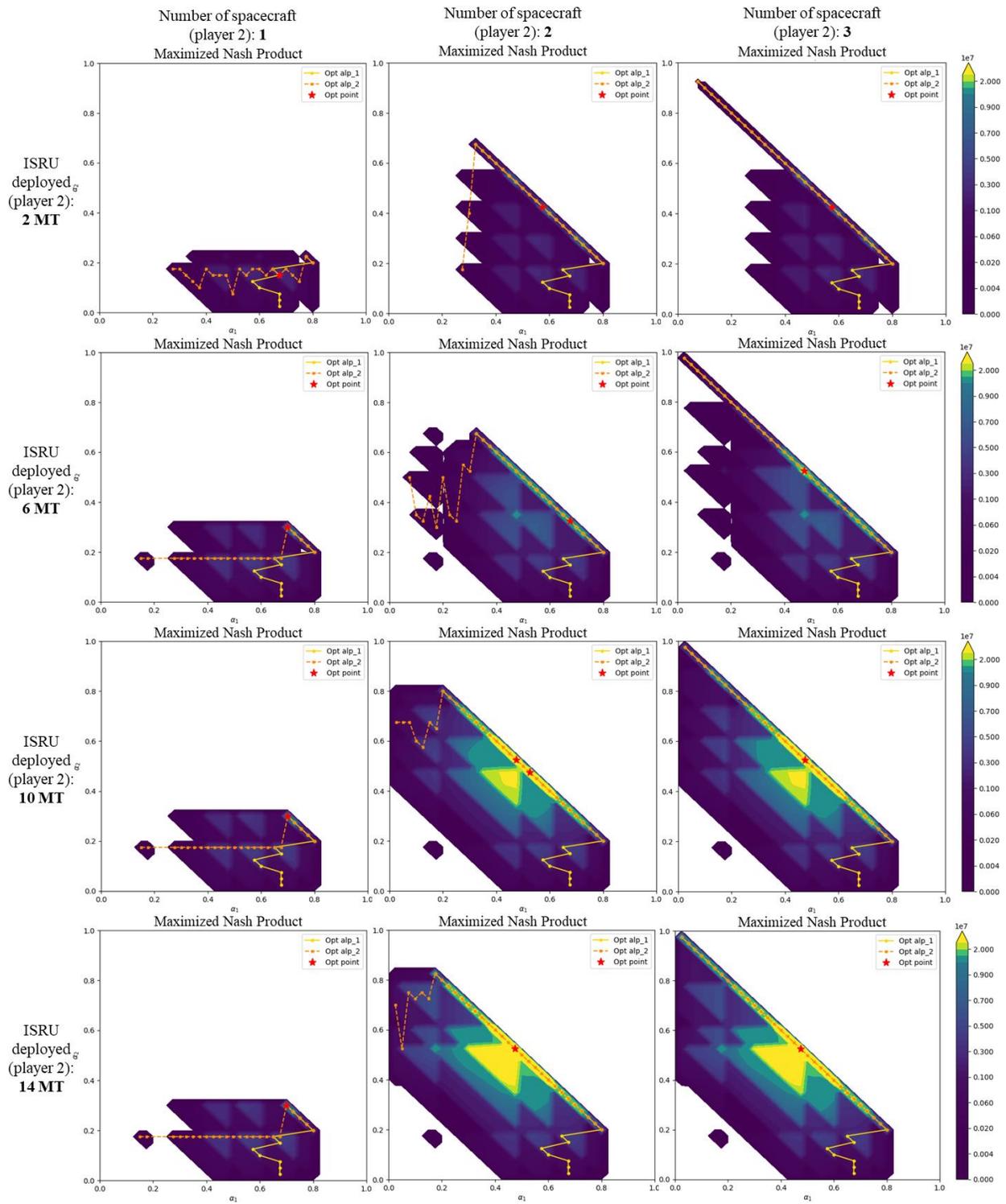

Fig. 11 Three-Player Design Space with Different Size of ISRU and Number of Spacecraft

## V. Conclusions

This paper proposes a space exploration architecture decision-making framework by extending state-of-the-art space logistics design methods from the game-theoretic perspective to enable the analysis in stimulating commercial participation in space infrastructure development and deployment. In the framework, a participation coefficient $\boldsymbol{\alpha}$ is defined to measure the deployment mission task assignment, and an incentive coefficient $\boldsymbol{\theta}$ is defined to measure the strength of the incentive. Based on the mission planning circumstances, there are three different scenarios for $\boldsymbol{\alpha}$ and $\boldsymbol{\theta}$ that correspond to the scenarios with the optimal design, fixed demand allocation, and fixed budget allocation. We analyze the incentive design properties under these scenarios and propose corresponding approaches to solve the problem.

A lunar exploration case study is performed to demonstrate the value of the incentive design framework. The analysis shows that by implementing an adequately defined incentive mechanism, the mission coordinator can complete its space mission at a lower total expense leveraging other commercial players' resources. Results show the value of space commercialization for future human space exploration. Moreover, sensitivity analyses also show the relationships between mission demands and the design space of the participation coefficient $\boldsymbol{\alpha}$, and between the size of the ISRU plant and the design space of the incentive coefficient $\boldsymbol{\theta}$. As an outcome of this case study, we derive an incentive design framework that can benefit both the mission coordinator and the commercial players from commercialization, leading to a mutually beneficial framework between the government and the industry. Furthermore, findings from the three-player mission analysis demonstrate the value of the proposed method in generating effective development strategies in a space commercialization environment.

Future research includes the consideration of uncertainties in the incentive design, including demand changes, spacecraft flight or rocket launch delay, and infrastructure performance uncertainties. Further implementations can also be explored to consider multiple commercial players with different technology strengths or take into account the long-term perspective of technology development by accepting temporary high mission costs.

## Acknowledgments

This paper is an extension to Ref. [25], which was presented at the 70[th] International Astronautical Congress, Washington DC in Oct. 2019.